\documentclass[a4paper,draft]{amsart}
\usepackage{amsmath}
\usepackage{amssymb}
\usepackage{bbm}
\usepackage{a4wide}

\parskip\smallskipamount

\def\<#1>{\langle#1\rangle}

\def\e{\mathrm{e}}

\newcommand\todo[1][.]{\edef\tmpa{.}\edef\tmpb{#1}%
  \ifx\tmpa\tmpb
    \typeout{To Be on page \thepage}\fbox{\bf To Be}
  \else
    \typeout{To Be on page \thepage: #1}\fbox{{\bf To Be:} #1}
  \fi
}

\begin{document}

 \author[Manuel Kauers and Doron Zeilberger]{Manuel Kauers\,$^1$ and Doron Zeilberger\,$^2$}

 \address{Manuel Kauers, Research Institute for Symbolic Computation, J. Kepler University Linz, Austria}
 \email{mkauers@risc.jku.at}

 \address{Doron Zeilberger, Mathematics Department, Rutgers University (New Brunswick), Piscataway, NJ, USA}
 \email{zeilberg@math.rutgers.edu}

 \thanks{$^1$ Partially supported by the Austrian FWF grant Y464-N18.}

 \thanks{$^2$ Partially supported by the USA National Science Foundation.}
 \date{ Nov. 21, 2010}
 \title{The Computational Challenge of Enumerating High-Dimensional Rook Walks}


 \maketitle

 \section{Introduction}

Consider a rook placed on the lower left corner $(0,0)$ of a chess board.
On how many paths can the rook reach the upper right corner $(n,n)$ if in a single step
it may move an arbitrary number of fields upwards or to the right (but not downwards
or to the left)? For rectangular chess boards of size $n\times m$, the number
of paths is given by the coefficient $a_{n,m}$ in the rational series expansion
\[
  \sum_{n,m=0}^\infty a_{n,m}x^n y^m=\frac1{1-\frac x{1-x}-\frac y{1-y}}.
\]
The case of square chess boards is consequently the \emph{diagonal series} of this rational function, which happens to be
\[
  \sum_{n=0}^\infty a_{n,n} x^n = \frac12+\frac{1-x}{2\sqrt{1-10x+9x^2}}.
\]
{}From here, all sorts of information about the numbers~$a_{n,n}$ can be easily extracted by 
means of computer algebra, for instance the initial terms
\[
  1,\ 2,\ 14,\ 106,\ 838,\ 6802,\dots\qquad\text{(A051708)},
\]
or recurrence equation 
\[
  (n+2)a_{n+2,n+2} - (10n+14)a_{n+1,n+1} + 9n\,a_{n,n} = 0\qquad(n\geq0),
\]
or the asymptotic formula $a_{n,n}\sim\sqrt{\!\frac2{\pi n}\!}\,3^{2n-1}$ ($n\to\infty$).
Computer algebra can also find the algebraic representation of the diagonal series given
the bivariate rational series as input, so there is altogether no need to do any calculation by hand.

At a marvelous meeting at Nankai University in August 2010 on the occasion of the second 
author's 60th birthday, Fr\'ed\'eric Chyzak reported that he and his colleagues had succeeded 
in doing the analogous computation for 3D~\cite{bostan10a}, i.e., they determined the number of paths
a rook can take on a 3D chessboard from $(0,0,0)$ to $(n,n,n)$ moving in each step
an arbitrary positive integer distance into one of the three directions, i.e., moving
either by $(i,0,0)$ or by $(0,i,0)$ or by $(0,0,i)$ for some positive integer~$i$.
Denoting now the number of this kind of walks by~$a_n$, they found the initial terms
\[
  1,\ 6,\ 222,\ 9918,\ 486924,\ 25267236,\dots\qquad\text{(A144045)},
\]
the recurrence equation
\begin{alignat*}1
&2(2 + n)(3 + n)^2(53 + 35n)a_{n+3}\\
&\quad-(2 + n)(43362 + 63493n + 30114n^2 + 4655n^3)a_{n+2}\\
&\qquad+(1 + n)(54864 + 100586n + 59889n^2 + 11305n^3)a_{n+1}\\
&\qquad\quad-192n^2(1 + n)(88 + 35n)a_n=0\qquad(n\geq0),
\end{alignat*}
and the asymptotic formula $a_n\sim\frac{9\sqrt3}{40\pi n}64^n$ ($n\to\infty$), all fully
rigorous, including so-called certificates which allow for an independent formal verification 
of the obtained results. 
While it is clear \emph{in theory} that computer algebra is able to obtain this information,
it is remarkable that it is possible to actually carry out these calculations \emph{in practice,} 
because the 3D case requires far more computational power than the 2D case.

If we don't insist on a fully rigorous formal verification, the diagonal recurrence can be obtained
with much less effort: it suffices to compute some 25 terms of the sequence and use \emph{automated guessing}
to find a recurrence which matches them. See \cite{bostan09,rubey07} for recent developments of this 
technique and the references given there for classical versions. 
For the present paper we applied this technique to empirically find recurrence equations for rook paths 
in dimensions greater than three, and we pose it as a challenge to provide rigorous certificates for
them. While at least for the very high dimensions this seems totally hopeless for now, we do expect that 
the coming years (or decades?) will see not only faster and bigger computers but also more advanced 
algorithms which can certify our claims within a reasonable amount of computing time. 
At least we intend to encourage progress in this direction. 
We see no other reason to ask for a certification. 
The question cannot be whether our claims are 
correct---the empirical evidence is way too strong to leave any reasonable doubt about that.
Nor can the question be whether there actually exist certificates for our 
claims---it is clear by theory that recurrences of diagonal sequences of multivariate rational series 
can always be certified. 
Nor can the question be whether a proof may provide some insight or 
understanding---certificates are usually just messy polynomials.
The interesting questions instead are: how big are the certificates, what is the computational cost 
for constructing them, and who will be the first to get the computation done. 

\section{An Alternative Route for Turning our Semi-rigorous proofs to Full-Fledged Rigorous Proofs}

We know {\it a priori} that there is a recurrence, this follows from general holonomic nonsense.
But by the work of  Moa Apagodu and Zeilberger \cite{az06} one can derive {\it a priori} upper bounds
for the promised recurrences. The recurrences for $d$-dimensional rook walks turn out
(empirically) to have order $d$ (for $2 \leq d \leq 12$). It is very possible that it  won't be too hard to
prove this sharp upper bound in general, or even a  weaker---but still realistic---one. This may enable one to
give a ``soft'' proof that the empirically ``guessed'' recurrences are indeed {\it rigorously}  proved.

If one would be able to find (realistic!) {\it a priori} bounds for the degrees of the coefficients as well,
then by plain {\it linear algebra} the ``guessed'' recurrences would be rigorously proved.

\section{A short interlude: Why is this problem So important?}

The harsh and/or sceptical critic may say: Who cares? Not that many people (or machines) play 12-dimensional chess,
and even the vast majority of the many people that {\tt do} play traditional 2D, $8 \times 8$ chess,
couldn't care less about the number of ways a rook can walk.

But {\it everyone} owes money, and usually to many creditors! The number of ways a rook can walk, in
the $d$-dimensional cubic lattice, from the origin to $[n, \dots, n]$ is also the number of ways
of repaying all your creditors if you currently owe $n$ dollars to each of $d$ different creditors, and
a single payment consists of paying any positive amount of dollars (up to the whole debt) to
any one of your creditors. {\it Now} this is a {\bf very} practical problem.

 \section{Fast Computation of Sufficiently Many Sequence Terms}

As the dimension increases, so does the order of the diagonal recurrence and the degree of the polynomials
appearing in it. The larger a recurrence is in terms of order and degree, the more sequence terms are needed 
to recover it from sequence data. For instance, in dimension $d=12$, we needed 1600 diagonal terms in order to find the 
recurrence. To modern guessing software (we used code written by the first author~\cite{kauers09a}), this is
still a moderate problem size. Much harder than guessing the recurrence is the computation of sufficiently 
many terms on the diagonals, which are needed as input for the guesser. The naive way is to start from the rational
function
\[
  \frac{p(x_1,\dots,x_d)}{q(x_1,\dots,x_d)}
  =\frac1{1-\frac{x_1}{1-x_1}-\frac{x_2}{1-x_2}-\cdots-\frac{x_d}{1-x_d}}
  =\sum_{n_1,\dots,n_d=0}^\infty a_{n_1,\dots,n_d}x_1^{n_1}\cdots x_d^{n_d}.
\]
Its denominator $q(x_1,\dots,x_d)$ gives rise to a multivariate linear recurrence with constant 
coefficients, which can be used to compute the $a_{n_1,\dots,n_d}$ recursively. For example, for $d=2$,
the rational function 
\[
  \frac1{1-\frac x{1-x}-\frac y{1-y}}=\frac{(x-1)(y-1)}{1-2x-2y+3xy}
\]
implies the recurrence
\[
  3a_{n+1,m+1}-2a_{n,m+1}-2a_{n+1,m}+a_{n,m}=0.
\]
Together with suitable boundary conditions, this allows the computation of $a_{n,m}$ for arbitrary
$n,m$, and hence for $a_{n,n}$ for arbitrary~$n$.

But this is very costly. In dimension~$d$, in order to compute the $n$th diagonal term, the recurrence
forces us to compute all terms $a_{n_1,\dots,n_d}$ with $0\leq n_1,\dots,n_d\leq n$, altogether 
more than $n^d$ terms. If $n=1000$, a computer won't mind doing this for $d=2$, but for $d=3$ 
it is already getting painful, and for $d>3$ either the memory requirements will exceed the computer's 
capacity or the runtime will exceed the user's patience. Or both. 
For $d\geq10$ the naive method will not even suffice for computing the first $n=10$ diagonal terms
within a reasonable amount of time.

Fortunately, there are more efficient recurrence equations. For arbitrary dimension~$d$, we have
\begin{alignat*}1
 n_d a(n_1,\dots,n_{d-2},n_{d-1},n_d)
 &=(n_{d-1}-1) a(n_1,\dots,n_{d-2},n_{d-1}-1,n_d-1)\\
 &\quad+(n_{d-1}+1) a(n_1,\dots,n_{d-2},n_{d-1}+1,n_d-1)\\
 &\qquad+(2-n_d) a(n_1,\dots,n_{d-2},n_{d-1},n_d-2)\\
 &\qquad\quad+(2 n_d-2 n_{d-1}-2) a(n_1,\dots,n_{d-2},n_{d-1},n_d-1).
\end{alignat*}
Note that any application of this recurrence leaves the indices $n_1,\dots,n_{d-2}$ fixed,
increases $n_{d-1}$ and decreases~$n_d$. This special form breaks the exponential complexity.
It can be shown that computing the first $n$ diagonal terms via this recurrence requires 
only $\mathrm{O}(n^2d^3)$ operations. For $d\leq7$, this method was efficient enough to produce
enough terms to obtain the recurrence for the diagonal. 

For $d\geq8$, an additional improvement was needed. Here instead of directly computing the terms
on the main diagonal, we first used the previous method for computing the terms of the bivariate 
auxiliary sequence
\[
  b_{n,m} := a_{n,\dots,n,m},
\]
up to $n,m\leq200$ or so. Then we used a multivariate guesser to discover some bivariate recurrences in
$n$ and $m$ for $b_{n,m}$ and used these to compute the diagonal terms $b_{n,n}=a_{n,\dots,n,n}$
for $n$ as far as needed.

\section{Recurrence Equations}

Most of the recurrences we found are too big to be reproduced here. We make them available online at
\begin{center}
  \verb|http://www.risc.jku.at/people/mkauers/walks/|.
\end{center}
Here we only give a table with some statistics on their order, the maximal degree of their 
polynomial coefficients, and the length of the longest integer appearing in them, measured
in decimal digits (dd).

\begin{center}
 \begin{tabular}{c|ccc|cl}
  dim & ord & deg & maxint & OEIS tag & comment \\\hline
  2 & 2 & 1 & 2\,dd & A051708 & easy\rule{0pt}{1.3em} \\
  3 & 3 & 4 & 6\,dd & A144045 & Chyzak et al.'s result\rule{0pt}{1.3em} \\
  4 & 4 & 9 & 12\,dd & A181749 &\rule{0pt}{1.3em}\\
  5 & 5 & 18 & 31\,dd & A181750 &\rule{0pt}{1.3em}\\
  6 & 6 & 31 & 51\,dd & A181751 &\rule{0pt}{1.3em}\\
  7 & 7 & 50 & 94\,dd & A181752 &\rule{0pt}{1.3em}\\
  8 & 8 & 75 & 149\,dd & A181754 &\rule{0pt}{1.3em}\\
  9 & 9 & 108 & 236\,dd & A181725 &\rule{0pt}{1.3em}\\
  10 & 10 & 149 & 306\,dd & A181726 &\rule{0pt}{1.3em}\\
  11 & 11 & 200 & 462\,dd & A181727 &\rule{0pt}{1.3em}\\
  12 & 12 & 261 & 609\,dd & A181728 &\rule{0pt}{1.3em}
 \end{tabular}
\end{center}

\section{Queens}

We have also computed recurrences for the analogous problem
of {\it Queen} walks, but so far we were only able to go up to dimension~5.
The relevant output can be found in the above-mentioned webpage.

\section{Higher Order Asymptotics}

The leading-term asymptotics for diagonals of rook walks has been derived by 
Martin Erickson, Suren Fernando, and Khang Tran~\cite{eft10}
using the powerful analytical method of Robin Pemantle and Mark C. Wilson \cite{pw04}.
It turns out to be:
$$
\sqrt{\alpha_d}(n\pi)^{(1-d)/2}((d+1)^{d})^{n} \quad \quad (n\to\infty)  \quad,
$$
where $\alpha_d$ is given by:
$$
\alpha_d = \frac{d^{d+2}}{(d+2)^{d-1}(d+1)^2 2^{d-1}} \quad .
$$
This result matches well with the numbers produced by the recurrences we discovered. 
For instance, for $d=12$ and $n=500000$ we find
\[
  \frac{\sqrt{\alpha_d}(n\pi)^{(1-d)/2}((d+1)^{d})^{n}}{\mathrm{A181728}(n)}=1.0000020411\dots
\]
And this not all. Thanks to the Maple package {\tt AsyRec} available from

\begin{center}
\verb|http://www.math.rutgers.edu/~zeilberg/tokhniot/AsyRec|
\end{center}

(see \cite{z08}) one can get, very easily, from the recurrences, higher-order asymptotics,
using the Birkhoff-Trjitznisky method.

The order-$10$ asymptotic formulas for  the sequences for $2 \leq d \leq 9$ can be gotten from 

\begin{center}
\verb|http://www.math.rutgers.edu/~zeilberg/tokhniot/oRookAsymptotics|
\end{center}

that is based on the input file

\begin{center}
\verb|http://www.math.rutgers.edu/~zeilberg/tokhniot/inRookAsymptotics|
\end{center}

that uses {\tt AsyRec} and of course, the recurrences obtained  by the first-named author's computer.

A cross check with $d=12$ and $n=500000$ now yields the very convincing quotient
\def\myBox{\raisebox{-.2ex}{$\Box$}}
\begin{alignat*}1
  &\frac{\sqrt{\alpha_d}(n\pi)^{(1-d)/2}((d+1)^{d})^{n}(1+\myBox\frac1 n+\myBox\frac 1{n^2}+\cdots+\myBox\frac 1{n^{10}})}{\mathrm{A181728}(n)}\\
  &\quad=0.999999999999999999999999999999999999999999999999999999999999999963446\dots
\end{alignat*}
where the $\Box$ symbol suppresses some explicit rational numbers which are too lengthy to be reproduced
here but which can be also found on the website above. 

And this is still not all. 
By looking at the output of {\tt AsyRec} for the sequences for {\it specific} $d$, it appears that 
we have the more refined asymptotic expression for the number of rook-walks from $[0^d]$ to
$[n^d]$ for fixed, but {\it arbitrary} (symbolic!) $d$
\begin{alignat*}1
  &\sqrt{\alpha_d}(n\pi)^{(1-d)/2}((d+1)^{d})^{n}\\
  &\quad\times\Bigl( \,\, 1- \frac{(d-1)(d+1)(d^3+6d^2+18d+12)}{12d(d+2)^3} \cdot \frac{1}{n} \\
  &\quad+{\frac { (d-1) (d+1)^2 ( {d}^{8}+11{d}^{7}+60{d}^{6}+168{d}^{5}-108{d}^{4}-564{d}^{3}-1632{d}^{2}-1584d-576)}{288{d}^{3} ( d+2 ) ^{6}}} \cdot \frac{1}{n^2}\\
  &\quad+ \mathrm{O}(\frac{1}{n^3})\,\, \Bigr)\quad.
\end{alignat*}

We leave the rigorous proof of this as another challenge to the reader.

\section{Fixed $n$, variable dimension}

Let $w_n(d)$ be the number of ways a rook can positively walk from $[0^d]$ to $[n^d]$.
So far, we fixed $d$ and let $n$ vary. But what if we fix $n$ and let $d$ vary?
Of course $w_0(d) \equiv 1$ and $w_1(d)=d!$, Sloane's A000142.
The sequence $w_2(d)$ is of more recent vintage, it is Bob Proctor's sequence A105749.
But a search on Nov. 19, 2010, did not find $w_3(d)$ in Sloane, or elsewhere.

The Maple package {\tt RookWalks} available from

\begin{center}
\verb|http://www.math.rutgers.edu/~zeilberg/tokhniot/RookWalks|
\end{center}

handles these sequences, and the webpage

\begin{center}
\verb|http://www.math.rutgers.edu/~zeilberg/tokhniot/oRookWalks|
\end{center}

lists the first 150 terms of $w_n(d)$ for $1 \leq n \leq 4$, as well as
guessed recurrences and implied asymptotics. The asymptotic
formulas for the individual $n$ (for $1 \leq n \leq 4$)
lead one to
conjecture that the leading asymptotics for $w_n(d)$  as $d \rightarrow \infty$ is
$$
\e^{n-1} \frac{(nd)!}{n!^d} \left ( \, 1 + \, \mathrm{O}(\frac{1}{d}) \right ) \quad .
$$
We leave the rigorous proof of this as yet another challenge to the reader.

\bibliographystyle{plain}
\bibliography{all}

\begin{thebibliography}{10}




\bibitem{az06} Moa Apagodu and Doron Zeilberger
\emph{Multi-Variable Zeilberger and Almkvist-Zeilberger Algorithms and the Sharpening of Wilf-Zeilberger Theory},
Adv. Appl. Math. {\bf 37}(2006), 139-152.

\bibitem{bostan09} Alin Bostan and Manuel Kauers,
\emph{Automatic classification of restricted lattice walks},
In Proceedings of FP-SACb09, pages 201b.

\bibitem{bostan10a} Alin Bostan, Frederic Chyzak, Mark van Hoeij, and Lucien Pech,
\emph{paper in preparation}.



\bibitem{eft10} Martin Erickson, Suren Fernando, and Khang Tran,
\emph{Enumerating Rook and Queen Paths},
Bulletin of the Institute of Combinatorics and its Applications. {\bf 60}(2010), 37-48.

\bibitem{rubey07} Waldemar Hebisch and Martin Rubey,
\emph{Extended Rate, more GFUN},  arXiv:math/0702086.

\bibitem{kauers09a} Manuel Kauers,
\emph{Guessing handbook},
Technical Report 09-07, RISC-Linz, 2009.

\bibitem{pw04} Robin Pemantle and Mark C. Wilson,
\emph{Asymptotics of multivariate sequences. II. multiple points of the singular variety},
Combinatorics, Probability, Computation {\bf 13} (2004), 735-761.



\bibitem{z08} Doron Zeilberger
\emph{AsyRec: A Maple package for Computing the Asymptotics of Solutions of Linear Recurrence Equations with Polynomial Coefficients},
Personal Journal of Shalsoh B. Ekhad and Doron Zeilberger, April 04, 2008.
\verb|http://www.math.rutgers.edu/~zeilberg/pj.html|.







\end{thebibliography}

\end{document}